\documentclass[12pt]{amsart}
\usepackage{}
\usepackage{cases}
\usepackage{graphicx}
\usepackage{amsmath}
\usepackage{amssymb}
\usepackage{amsfonts}
\usepackage[mathscr]{eucal}
\makeatletter
\def\numberwithin#1#2{\@ifundefined{c@#1}{\@nocnterrr}{
  \@ifundefined{c@#2}{\@nocnterr}{
  \@addtoreset{#1}{#2}
  \toks@\expandafter\expandafter\expandafter{\csname the#1\endcsname}
  \expandafter\xdef\csname the#1\endcsname
    {\expandafter\noexpand\csname the#2\endcsname
     .\the\toks@}}}}
\makeatother
\numberwithin{equation}{section}

\newtheorem*{problem*}{Problem 1}

\begin{document}
\title[ Simon's conjecture]
{Curvature bounds for the spectrum of closed Einstein spaces  and Simon conjecture}

\thanks{The authors are  supported by the National Science Foundation of China
(No.12161092).}
\author{ShanLin Guan}
\address {Department of Mathematics, Yunnan Normal University, Kunming, 650500, P. R. of China; School of Mathematical Sciences, FuDan University, 200433, Shanghai, P. R. of China}
\email{22110180016@m.fudan.edu.cn}
\author {Zhen Guo}
\address{Department of Mathematics, Kunming University, Kunming, 650214; Department of Mathematics, Yunnan Normal University, Kunming, 650500, P. R. of China}
\email{gzh2001y@yahoo.com}

\subjclass[2020]{53C24, 53C25, 53C21, 53C80}
\date{\today}

\keywords{Einstein space, Spectrum, Simon's conjecture}
\begin{abstract}
Let $(M^{n}, g)$ be a closed connected Einstein space, $n=dim M ,$ and $\kappa_{0} $ be the lower bound of the sectional
curvature. In this paper, we prove Udo Simon's conjecture: on closed Einstein spaces, $n\geq  3,$ there is no eigenvalue $\lambda$
such that
$$n\kappa_{0} < \lambda < 2(n + 1)\kappa_{0},$$
and both bounds are the best possible. Furthermore, we develop Simon's conjecture to the next gap of eigenvalue $\lambda:$ on closed Einstein spaces, there is no $\lambda$ such that
$$
 2(n + 1)\kappa_{0}< \lambda < 2(n+2)\kappa_{0}, 
$$
and both bounds are the best possible.
\end{abstract}

\maketitle

%%%%%%%%%%%%%%%%%%%%%%%%%%%%
%%%%%%%%%%%%%%%%%%%%%%%%%%%%

\section{Introduction}
 Purpose of this paper is to investigate the curvature bounds for the spectrum of closed Einstein spaces. The classical 
Lichnerowicz-Obata's theorem states that if $(M, g)$ is a $n$-dimensional closed Riemannian manifold with metric $g,$
where the Ricci curvature satisfies $Ric \geq  (n- 1)Kg$ for a constant $K,$ then
the first nonzero eigenvalue $\lambda_{1}$ of the Laplace operator satisfies $\lambda_{1}\geq nK$; and
the equality holds if and only if $M $ is isometric
to the sphere $\mathbb{S}^{n}(K)$. 
In the case of Einstein manifold, Lichnerowicz-Obata's  theorem was states as follows:
\par\noindent
\textbf{Theorem A} ( Lichnerowicz-Obata\cite{{Lich}, {O}} ) \emph{Let $M^{n}(n\geq 3)$ be a closed Einstein manifold with  scalar curvature $R$ and let $\lambda$ be the first eigenvalue of the Laplacian. Then 
$$\lambda_{1}\geq \frac{R}{n-1},$$
and the equality case holds if and only if $M^{n}$ is isometrically diffeomorphic to a sphere.}
\par\noindent
Tanno proved the following.
\par\noindent
\textbf{Theorem B} (\cite{Tan2}) \emph{Let $M^{n}(n\geq 2)$ be a closed Einstein manifold with  scalar curvature $R$ and let $\kappa_{0}$ be the lower bound of the sectional
curvature. Then there is no eigenvalue $\lambda$ of the Laplacian such that
$$\frac{R}{n-1}< \lambda\leq 2(n-1)\kappa_{0}+\frac{R}{n(n-1)},$$
}
Udo Simon proved following result. 
\par\noindent
\textbf{Theorem C} (\cite{Sim})  \emph{Let $(M^{n}, g)$ be a closed connected Einstein space, $n=dim M \geq 2$ (with constant scalar curvature $R$). Let $\kappa_{0} $ be the lower bound of the sectional
curvature. Then either $ (M, g) $ is isometrically diffeomorphic to a sphere and
the first nonzero eigenvalue $\lambda_{1}$  of the Laplacian fulfills
$$\lambda_{1} = \frac{ R }{n-1}= n\kappa_{0}$$
or each eigenvalue $\lambda$ of the Laplacian satisfies the inequality
$$\lambda > 2n\kappa_{0}.$$ }

 As for a sphere of constant sectional curvature $K,$  the first nonzero eigenvalue is given by $\lambda_{1} = nK$, the second by $\lambda_{2} = 2(n + 1)K.$  
 Basic on Theorem C, the examples in \cite{CW} and the value of $\lambda_{2}$ on spheres, Udo Simon proposed  the following.
\par\noindent
\textbf{ Conjecture }(\cite{Sim})\emph{ On closed Einstein spaces, $n\geq  3,$ there is no eigenvalue $\lambda$
such that
$$n\kappa_{0} < \lambda < 2(n + 1)\kappa_{0}.$$
Both bounds are the best possible.}
\par\noindent 
For resolving this conjecture,  K. Benko et al. used new estimates and obtained  the following.
\par\noindent
\textbf{Theorem D} ( cf. Theorem 3.1 in \cite{BKSS} ) \emph{On closed Einstein spaces with dimension $n\geq  3,$  $\lambda>n\kappa_{0}$ implies 
$$\lambda\geq 2(n+1)\kappa_{0}-2(\frac{R}{n(n-1)} -\kappa_{0}).
$$
}
\par\noindent
As far as we know, since then there is no progress on study of this conjecture. it is still open. 
In this paper, we confirm this conjecture. Explicity, we prove the following.
\par\noindent
\textbf{ Theorem 1.1} \emph{Let $(M^{n}, g)$ be a closed connected Einstein space, $n=dim M \geq 3,$ then  $\lambda >n\kappa_{0}$ implies 
$$\lambda\geq 2(n+1)\kappa_{0},$$
and the equality holds only if the universal covering space of $(M^{n}, g)$ is isometrically diffeomorphic to a sphere $\mathbb{S}^{n}(\kappa_{0}).$ }
\par
In fact, theorem A shows that the first gap of $\lambda$ is interval $(0, n\kappa_{0}),$ Theorem 1.1 gives the second gap $(n\kappa_{0}, 2(n+1)\kappa_{0})$ on a closed Einstein space with positive $\kappa_{0}.$ Furthermore, we investigate the next gap of the eigenvalue $\lambda$  and get the following.
\par\noindent
\textbf{ Theorem 1.2} \emph{
		Let $(M^{n}, g)$ be a closed connected Einstein space and $\kappa_{0} $ be the lower bound of the sectional
		curvature, $n=dim M \geq 3$. Then $\lambda > 2(n+1) \kappa_{0}$ implies $$\lambda \ge 2(n+2)\kappa_{0}.$$ 
	    and the equalitiy holds implies the universal covering space of $M^{n}$ is  $\mathbb{CP}^{m}(c)$ with  constant holomorphic sectional curvature $c=4\kappa_{0}$, where $n=2m$.}
\par
The key step of the proof of the Simon's conjecture is to prove Lemma 3.1. Where we construct the minimal norm tensor of the third covariant derivatives of the eigenfunction to estimate the necessary inequality. \section{notations and formulae  }
Let $(M^{n}, g)$ be a connected Riemannian manifold of class $C^{\infty}.$  Let $\nabla$ be the covariant differentiation induced by $g;$ $R_{ijkl}, R_{ij}$  be the components of the Riemannian curvature and the Ricci curvature in local orthonormal basic $\{e_{i}, 1\leq i\leq n\},$  respectively, and $R$ be the scalar curvature. Let $f$ be a smooth function in $M^{n},$ $f_{i}, f_{i,j}$ and $f_{i,jk}$ be the components of the first, second and third order coderivatives of $f$ in local orthonormal basis $\{e_{i}, 1\leq i\leq n\},$ respectively;  $\Delta f,$ $\nabla f$ and $Hess(f)$ be the Laplacian, gradient  and Hessian of $f$, respectively.  We use following notations: 
$$
\|\nabla f\|^{2}=\sum_{i}f_{i}^{2},\ \
\|Hess(f)\|^{2}=\sum_{i,j}f_{i,j}^{2}, \ \
\|\nabla^{3} f\|^{2}=\sum_{i,j,k}f_{i,jk}^{2}.
\eqno(2.1)$$
On local orthonormal basis, the Ricci identity  can be written as follows:
$$
f_{i,jk}-f_{i,kj}=f_{l}R_{lijk}
\eqno(2.2)
$$
Let $(M^{n}, g)$ be an Einstein space, denote the Weyl conformal curvature tensor by $W,$ then we have the following formula:
$$
W_{ijkl}=R_{ijkl}-\frac{R}{n(n-1)}(\delta_{ik}\delta_{jl}-\delta_{il}\delta_{jk}),
\eqno(2.3)
$$
where $R_{ij}=\sum_{k}R_{ikjk}$ and $R=\sum_{i}R_{ii}.$
\section{Some Lemmas}
\par\noindent
\textbf{ Lemma 3.1 } \emph{
		Let $M^{n}$ be a $n-$dimensional  Einstein space and $f$ be the eigenfunction corresponding to eigenvalue $\lambda,$ then we have the following inequality.
$$
\|\nabla^{3} f\|^{2} \geq \left[ \frac{\lambda^{2}}{n} + \frac{2(n-1)}{n(n+2)} (\frac{R}{n-1} - \lambda)^2\right] \|\nabla f\|^{2} + \frac{1}{3} f_{i}W_{ipqr}W_{jpqr}f_{j},
\eqno(3.1)
$$  
and the equality holds if and only if
$$
\begin{array}{lc}
f_{i,jk} \displaystyle{= -\frac{1}{n+2}[\frac{2R}{(n-1)n} + \lambda]f_{k} \delta_{ij} - \frac{1}{n+2}(\lambda - \frac{R}{n-1})(f_{j}\delta_{ik} + f_{i}\delta_{jk}) }
\\
\displaystyle{+ \frac{1}{3}(f_{l}W_{lijk} - f_{l}W_{ljki})}.
\end{array}
\eqno(3.2)$$
}
\par\noindent
\emph{Proof.} Defining  tensor $A,$ $B$ and $F$ by
$$
A_{ijk}=\frac{1}{n+2}[\frac{2R}{(n-1)n} + \lambda]f_{k} \delta_{ij} + \frac{1}{n+2}(\lambda - \frac{R}{n-1})(f_{j}\delta_{ik} + f_{i}\delta_{jk}) ,
$$
$$
B_{ijk}=-\frac{1}{3}(f_{l}W_{lijk} - f_{l}W_{ljki}),
$$
$$
F_{ijk}=f_{i,jk} +A_{ijk}
+ B_{ijk}.
$$
It is easy to check that both tensor $\nabla^{3}f+A$ and $B$ are trace-free, and so $F$ is trace-free tensor. This implies that both $F$ and $B$ are orthogonal to $A,$ i.e. $<F,A>=<B,A>=0.$ From
$$
F-A=\nabla^{3}f+B,
$$
we have
$$\|F\|^{2}=\|\nabla^{3}f\|^{2}+2<\nabla^{3}f, B>+\|B\|^{2}-\|A\|^{2},
\eqno(3.3)$$
$$
<\nabla^{3}f, A>=-\|A\|^{2}.
\eqno(3.4)$$
We have following inequality:
$$
\|\nabla^{3}f\|^{2}\geq \|A\|^{2}-2<\nabla^{3}f, B>-\|B\|^{2}.
\eqno(3.5)
$$
where the equality holds if and only if $F=0.$ Using Ricci identity and (2.3),  we have the following computation.   
\[
\begin{aligned}
 <\nabla^{3}f, B> & =-\frac{1}{3}(f_{p}W_{pijk}f_{i,jk}-f_{p}W_{pjki}f_{j,ik})\\
  & =-\frac{1}{3}[f_{p}W_{pijk}f_{i,jk}-f_{p}W_{pjki}(f_{j,ki}+f_{q}R_{qjik})] \\
  & = -\frac{1}{3}f_{p}W_{pijk}R_{qijk}f_{q}\\
  & =-\frac{1}{3}f_{p}W_{pijk}W_{qijk}f_{q}. \end{aligned}
  \eqno(3.6) 
 \] 
 Noting that $W$ has of properties: $W_{pijk}=-W_{pikj},$ and $W_{pijk}+W_{pjki}+W_{pkij}=0,$ we have
 $$
  W_{pijk}W_{qjik}=\frac{1}{2}W_{pijk}W_{qijk}.
  \eqno(3.7)
  $$
  By using (3.7), we get
  $$
  \|B\|^{2}=\frac{1}{3}f_{p}W_{pijk}W_{qijk} f_{q}.
  \eqno(3.8)
   $$
   From (3.4) we have 
 $$
   \|A\|^{2}=-\sum_{ijk}f_{i,jk}A_{ijk}
   =\left[ \frac{\lambda^{2}}{n} + \frac{2(n-1)}{n(n+2)} (\frac{R}{n-1} - \lambda)^2\right] \|\nabla f\|^{2}.
 \eqno(3.9)
 $$
  Putting (3.6), (3.7) and (3.9) into (3.5),  we have (3.1). As the equality of (3.5) holds if and only if  $F=0,$ the equality holds of (3.1) if and only if $f$ satisfies (3.2). This completes the proof of Lemma 3.1.
  \par\noindent
\textbf{ Lemma 3.2} \emph{
Let $M^{n}$ be a $n-$dimensional closed Einstein space, then we have the following equality.
 \[
\begin{aligned}
\int_{M^{n}} \frac{1}{3}f_{i}W_{ipqr}W_{jpqr}f_{j} = \frac{2}{3} \int_{M^{n}} \sum_{i < j} [\frac{R}{n(n-1)} -\kappa_{ij}](h_{i} - h_{j})^{2}. 
\end{aligned}
\eqno(3.10)	\]	
 where $h_{i}$ is the eigenvalues of $Hess(f),$ $\{e_{i}\}$ is corresponding orthonormal basis and
$\kappa_{ij}$ is the sectional curvature of the plane $\{e_{i}, e_{j}\}_{i\neq j}.$   } 
  
\par\noindent
\emph{Proof.} 
By using Ricci identity we  have
$$
f_{q}W_{qijk}W_{pijk}=2f_{i,jk}W_{pijk}. \eqno(3.11)
$$
Noting that $\sum_{k}W_{pijk,k}=0 $ holds on an Einstein space, Using (3.11), we have
    	\[
	\begin{aligned}
    			\frac{1}{3} \int_{M^{n}}  f_{l}W_{lijk}W_{pijk}f_{p} & = \frac{2}{3} \int_{M^{n}}  f_{i,jk}W_{pijk}f_{p}  \\ &  = - \frac{2}{3}\int_{M^{n}}  f_{i,j}W_{pijk}f_{p,k}. 
			\end{aligned}
		\eqno (3.12)	\]
At a point we can take $e_{i}$ such that $f_{i,j}=h_{i}\delta_{ij}.$ Noting that $W$ is trace free, we have
$$
		\begin{aligned}		
		f_{i,j}W_{pijk}f_{p,k}=-W_{ijij}h_{i}h_{j} & =\sum_{i<j}W_{ijij}(h_{i}-h_{j})^{2}\\
	&=\sum_{i<j}[\kappa_{ij}-\frac{R}{n(n-1)}](h_{i}-h_{j})^{2}.
		\end{aligned}
		\eqno(3.13)
$$
This completes the proof of Lemma 3.2.
\par\noindent
\textbf{ Lemma 3.3} \emph{
		Let $M^{n}$ be a $n-$dimensional closed Einstein space, then we have the following equality.}
\[
\int_{M^{n}} \sum_{i < j} (h_{i} - h_{j})^2= [(n-1)\lambda - R] \int_{M^{n}} \|\nabla f \|^{2}. \eqno(3.14)
\]
\par\noindent
\emph{Proof.} 
On a closed Einstein space, we have following identity:
$$
\int_{M^{n}}[n\|Hess(f)\|^{2}-(\Delta f)^{2}]=
[(n-1)\lambda - R] \int_{M^{n}} \|\nabla f \|^{2}.
\eqno(3.15)
$$	
Noting 
$$
n\|Hess(f)\|^{2}-(\Delta f)^{2}=\sum_{i < j} (h_{i} - h_{j})^2,
\eqno(3.16)
$$
we have (3.14). 
\par\noindent
\textbf{ Lemma 3.4} (\cite{O}, \cite{Tan}) \emph{
 Let (M, g) be a complete and simply connected Riemannian
manifold. In order for (M, g) to admit a non-constant function f satisfying 
$$
f_{i,jk}+2c\delta_{ij}f_{k}+c(\delta_{ki}f_{j}+\delta_{jk}f_{i})=0. \eqno(3.17)
$$
for some positive constant $c$, it is necessary and sufficient that (M, g) is isometric to an Euclidean sphere $S^{n}(c)$.}
\par
In recently, some profound rigidity theorems on Einstein space have been proved (cf., for instance, \cite{Brd}, \cite{MW}, \cite{Yang}, \cite{Har}). Here, for the purpose to prove Theorem 1.2, we  use a rigidity result given by Xu and Gu (\cite{XG}) to show  the following. 
\par\noindent
\textbf{Lemma 3.5}\emph{
	Let $M^{n}(n \ge 3)$ be a closed Einstein manifold with lower bound of sectional curvature $\kappa_{0}$. If the scalar curvature satisfies
	$$n(n-1) \kappa_{0} \le R \le n(n+2) \kappa_{0}, $$
	then either $R = n(n-1) \kappa_{0} $ and the universal covering space of $M^{n}$ is isometric to $\mathbb{S}^{n}(\kappa_{0})$ with constant sectional curvature $\kappa_{0},$ or $R = n(n+2) \kappa_{0}$ and the universal covering space of $M^{n}$ is isometric to $\mathbb{CP}^{m}(4\kappa_{0})$ with constant holomorphic sectional curvature $4\kappa_{0}$, where $n=2m$.}

\begin{proof}
	If $n = 3$, then $M^{n}$ has of constant sectional curvature and so the universal covering space of $M^{n}$ is  $\mathbb{S}^{3}(\kappa_{0})$ with constant sectional curvature $\kappa_{0}$. Now we consider the case $n \ge 4$. From the assumption 
	$n(n-1) \kappa_{0} \le R \le n(n+2) \kappa_{0}, $
	we have 
	$$
	\begin{aligned}
		\frac{R}{n(n-1)} \ge \kappa_{0} \ge  \eta_{n}\frac{R}{n(n-1)}, 
\end{aligned}
$$
where $\eta_{n}=1-3/(n+2).$	
 From the rigidity theorem (Theorem 1.1 in \cite{XG}) given by Xu and Gu we know the universal covering space of $M^n$ is isometric to either the $n$-sphere $\mathbb{S}^{n}(\kappa_{0})$ with constant sectional curvature $\kappa_{0}$ or the
	complex projective space $\mathbb{CP}^{m}(c)$ with constant holomorphic sectional curvature $c$, where 
	$$ n= 2m, c = \frac{4R}{n(n+2)}. $$
	For the case of $\mathbb{S}^{n}(\kappa_{0})$,  we have $R=n(n-1)\kappa_{0}$. For the case of $\mathbb{CP}^{m}(c)$ with constant holomorphic sectional curvature $c$, since  
$ \kappa_{0}=c/4,$ we have $ R = n(n+2)\kappa_{0} .$ This proves Lemma 3.5.
\end{proof}
	
\section{Proofs of Theorem 1.1 and Theorem 1.2}
\par
\emph{Proof of Theorem 1.1} \ \ Using Lemma 3.1 and Lemma 3.2, we have
$$
\begin{aligned}
\int_{M^{n}} \|\nabla^{3}f \|^{2} &\geq \int_{M^{n}}[\frac{\lambda^{2}}{n} +
\frac{2(n-1)}{n(n+2)}(\lambda-\frac{R}{n-1})^{2}]\|\nabla f\|^{2}\\
& +\frac{2}{3}\int_{M^{n}}\sum_{i<j}[\frac{R}{n(n-1)}-\kappa_{ij}](h_{i}-h_{j})^{2}.
\end{aligned}
\eqno(4.1)
$$
The next integral formula was given by Udo Simon (\cite{Sim},(2.4c))
$$
\int_{M^{n}}[\|\nabla^{3}f\|^{2}-\frac{\lambda^{2}}{n}\|\nabla f\|^{2}+\sum_{i<j}(2\kappa_{ij}-\frac{\lambda}{n})(h_{i}-h_{j})^{2}
]=0.\eqno(4.2)
$$
 From (4.1) and (4.2) we have
 $$
 0\geq \int_{M^{n}}[\frac{2(n-1)}{n(n+2)}(\lambda-\frac{R}{n-1})^{2}\|\nabla f\|^{2}
 +\sum_{i<j}[\frac{2R}{3n(n-1)}+\frac{4}{3}\kappa_{ij}-\frac{\lambda}{n}](h_{i}-h_{j})^{2}].
 \eqno(4.3)
  $$	
Using Lemma 3.3 and (4.3) we get
$$
\begin{aligned}
 0 & \geq \int_{M^{n}}\sum_{i<j}[-\lambda +\frac{2}{3n}R+\frac{4(n+2)}{3}\kappa_{ij}](h_{i}-h_{j})^{2}\\
& =\int_{M}[-\lambda+2(n+1)\kappa_{0}]\sum_{i<j}(h_{i}-h_{j})^{2}\\
&+\int_{M}\frac{2(n-1)}{3}[\frac{R}{n(n-1)}-\kappa_{0}]\sum_{i<j}(h_{i}-h_{j})^{2}\\
&+\int_{M}\frac{4(n+2)}{3}\sum_{i<j}(\kappa_{ij}-\kappa_{0})(h_{i}-h_{j})^{2}\\& \geq \int_{M}[-\lambda+2(n+1)\kappa_{0}]\sum_{i<j}(h_{i}-h_{j})^{2}.
\end{aligned}
\eqno(4.4)
$$
If $\sum_{i<j}(h_{i}-h_{j})^{2}=0,$ then $f_{i,j}=-\lambda f\delta_{ij}.$ So,  Lemma 3.3 and Theorem A show that $M^{n}$ is isometrically diffeomorphic to a sphere with curvature $\kappa_{0},$ the first eigenvalue $\lambda=n\kappa_{0}.$ Otherwise, we have $\lambda\geq 2(n+1)\kappa_{0}.$ From (4.4) we see that $\lambda= 2(n+1)\kappa_{0}$ implies $\frac{R}{n(n-1)}=\kappa_{0},$ $\sum_{i<j}(\kappa_{ij}-\kappa_{0})(h_{i}-h_{j})^{2}=0$
and the equality holds of (3.1). From Lemma 3.2 we have  $\sum_{l}f_{l}W_{lijk}=0$ and $f$ satisfies the equation (3.2),  which means $f$ satisfies 
$$
f_{i,jk} \displaystyle{= -2\kappa_{0}f_{k} \delta_{ij} - \kappa_{0}(f_{j}\delta_{ik} + f_{i}\delta_{jk}) }.
\eqno(4.5)
$$
Lemma 3.4 shows that the universal space of 
 $(M^{n}, g)$ is isometrically diffeomorphic to a sphere with curvature $\kappa_{0}$. This completes the proof of the Theorem 1.1.
 
\begin{proof} [Proof of Theorem 1.2] From Lemma 3.5 we know that either $R = n(n-1) \kappa_{0} $ or $R \ge n(n+2) \kappa_{0}.$ For the first case,  all sectional curvatures must be the same i.e. $\kappa_{ij} = \kappa_{0}$ and so the  universal covering space of $M^{n}$ is isometric to  the $n$-sphere $\mathbb{S}^{n}(\kappa_{0}).$  For the second case, making use of  the first inequality in (4.4), we have  
$$
\begin{aligned}
0 & \geq \int_{M^{n}}\sum_{i<j}[-\lambda +\frac{2}{3n}R+\frac{4(n+2)}{3}\kappa_{ij}](h_{i}-h_{j})^{2}\\
& \geq \int_{M^{n}}\sum_{i<j}[-\lambda +2(n+2)\kappa_{0}](h_{i}-h_{j})^{2}
\end{aligned}
		\eqno(4.5)
$$
 If $\sum_{i<j}(h_{i}-h_{j})^{2}=0,$ Then, as argument as that in the proof of Theorem 1.1 we know $f$ is the first eigenfunction and $\lambda$ is the first eigenvalue $n\kappa_{0}$ of sphere $\mathbb{S}^{n}.$ Hence, the condition $\lambda>2(n+1)\kappa_{0}$ of Theorem 1.2 and (4.5) shows
 $$
 \lambda\geq 2(n+2)\kappa_{0},
 $$ 
 and the equality holds if and only if $R=n(n+2)\kappa_{0}.$ From Lemma 3.5 we know that the universal covering space of $M^{n}$ is isometric to $\mathbb{CP}^{m}(4\kappa_{0})$ with constant holomorphic sectional curvature $4\kappa_{0}$, where $n=2m$. This completes the proof of Theorem 1.2.
 \end{proof}
 %%%%%%%%%%%%%%%%%%%%%%%%%%%%
%%%%%%%%%%%%%%%%%%%%%%%%%%%%
%%% \bibliographystyle{amsplain}

\end{document}